\documentclass[twoside, 11pt]{article}
\usepackage{amsmath} 
\usepackage[e]{esvect}
\usepackage{graphicx}
\usepackage{epstopdf}
\usepackage[adobe-utopia]{mathdesign}
\usepackage[T1]{fontenc}
\usepackage[english]{babel}
\usepackage{color}

\usepackage{hyperref}

\title{\textsc{There are at most $2^{d+1}-2$ neighbourly simplices in dimension $d$}}
\date{}
\author{Andrzej P. Kisielewicz \& Krzysztof Przes{\l}awski}
{

\pagestyle{myheadings}
\markboth{\hfill\sc neighbourly simplices \hfill}{\hfill\sc Kisielewicz \,\, Przes{\l}awski \hfill}
\newtheorem{lemat}{\sc Lemma}
\newtheorem{tw}{\sc Theorem}

\newtheorem{df}{\sc Definition}

\newtheorem{uw}{\sc Remark}
\newtheorem{uwi}[uw]{\sc Remarks}

\newtheorem{nap}{\sc Example}

\newtheorem{nps}[nap]{\sc Examples}

\def\ka #1{\mathscr{#1}}
\def\kal #1 #2{\mathscr{#1}^{#2}}
\def\proof{\noindent \textit{Proof.\,\,\,}}
\def\skok{\vskip 0.1in\noindent}

\def\er{\mathbb{R}}

\def\Aut #1 #2{\operatorname{Aut}^{#1} (#2)}

\def\prop #1{\operatorname{prop}\, #1}

\def\bred #1 {\colorbox{red}{ #1}}
\def\red #1 {{\color{red} #1 }} 

\graphicspath{{figpdfkolor/}}
\DeclareGraphicsExtensions{.pdf,.png} 
\begin{document}
\maketitle

\begin{abstract}
A combinatorial theorem on families of disjoint sub-boxes of a discrete cube, which implies that there are at most $2^{d+1}-2$ neighbourly simplices in $\er^d$, is presented.
\end{abstract}

\section{Introduction}
A family of $d$-dimensional simplices in $\er^d$ is \textit{neighbourly} if the intersection of every two members is $(d-1)$-dimensional. It has been repeatedly conjectured that the maximum cardinality of such a family is $c_d=2^{d}$ (see \cite{Zaks1} for further references). The conjecture is verified up to dimension $3$ only. F. Bagemihl \cite{Bagemihl} proved that $8\le c_3\le 17$. V. Baston \cite{Baston} proved $c_3\le 9$. The final step $c_3=8$ was made by J. Zaks \cite{Zaks2}. The same author \cite{Zaks1} showed by a clever construction that $c_d\ge 2^d$. It was M. Perles \cite{Perles} who had found $c_d\le 2^{d+1}$. A slightly better estimate $c_d\le 2^{d+1}-1$ is shown in \cite[Chapter 14]{AZ}. (This chapter together with a recent post \cite{Kalai} on G. Kalai's blog are a great introduction to the subject of neighbourly families.) One of our goals is to prove that $c_d\le 2^{d+1}-2$. Basically, we shall follow Baston's approach with the combinatorial flavour added by Perles.

Let $\ka F$ be a neighbourly family in $\er^d$, $d\ge 2$. Let us arrange all the hyperplanes spanned by the facets of simplices belonging to $\ka F$ into a sequence $H_1, \ldots, H_n$. Each $H_i$ splits $\er^d$ into two halfspaces. Let us call them $H_i^0, H_i^1$. For every $\sigma\in \ka F$, let us define a unique word $v=v_1\cdots v_n$ of length $n$ over the alphabet $\{0,1,*\}$ as follows
$$
v_i=
\begin{cases}
0 & \text{if $H_i$ is spanned by a facet of $\sigma$ and $\sigma\subset H_i^0$},\\
1 & \text{if $H_i$ is spanned by a facet of $\sigma$ and $\sigma\subset H_i^1$},\\
* & \text{otherwise.}
\end{cases}
$$
Let $V$ be the set of all just defined words $v$. As is easily seen, $V$ satisfies the assumptions of our Theorem \ref{glowne} with $k=d+1$. Therefore, $|\ka F|=|V|\le 2^{d+1}-2$, as expected. 
\section{Main result}

A key observation concerns boxes contained in $\{0,1\}^n$. It is a particular case of \cite[Lemma 8.1]{KisPrz}.

Let $B=B_1\times\cdots \times B_n$ be a box contained in $\{0,1\}^n$. Let
$$
\prop B=\{i \colon B_i\neq \{0,1\}\}.
$$
Two boxes $B$ and $C$ contained in $\{0,1\}^n$ are said to be \textit{equivalent} if $\prop B =\prop C$. We shall need a kind of order relation: $A\preceq B$ if $\prop A\supseteq \prop B$. 

\begin{lemat}
\label{minbox}
Given a family $\ka B$ of disjoint boxes contained in $\{0,1\}^n$. Suppose $A\in \ka B$ is minimal with respect to $\preceq$. Let $[A]$ consists of all members $\ka B$ that are equivalent to $A$. If $\ka B$ is a tiling of $\{0,1\}^n$, then 
$$
|[A]_e|=|[A]_o|,
$$
where 
$[A]_e=\{B\in [A]\colon |\{i\in \prop A: A_i\neq B_i\}|\equiv 0\pmod 2\}$ and $[A]_o=\{B\in [A]\colon |\{i\in \prop A: A_i\neq B_i\}|\equiv 1\pmod 2\}$.

In particular, there is $B\in [A]$ such that the set $\{i\in \prop A: A_i\neq B_i\}$ is of odd cardinality. 
\end{lemat}

\proof
Let us define a sequence of functions $f_i\colon \{0,1\}\mapsto \{-1,1\}$, $i=1,\ldots, n$, as follows:
$$
f_i(x)=\begin{cases} 
(-1)^x & \text{for $i \in \prop A$,}\\

1 & \text{for $i \notin \prop A$}.
\end{cases}
$$
Let $f=f_1\otimes\cdots\otimes f_n$; that is $f(x_1,\ldots, x_n)= f_1(x_1)\cdots f_n(x_n)$. For every $X\subseteq \{0,1\}^n$, let us set $\sum_X f=\sum_{x\in X} f(x)$. It is easy to calculate that $\sum_Af\neq 0$. (Namely, $\sum_Af=(-1)^s2^{d-|\prop A|}$, where $s$ is the cardinality of the set $\{i \colon A_i=\{1\}\}$). Moreover, by the minimality of $A$ and the definition of $f$ we have,
$$ 
\sum_B f=
\begin{cases}
\sum_A f & \text{if $B\in [A]_e$,}\\
-\sum_A f& \text{if $B\in [A]_o$,}\\
0 & \text{if $B\in \ka B \setminus [A]$.}
\end{cases}
$$
Since also $\sum_{\{0,1\}^n} f=0$ and $\ka B$ is a partition of $\{0,1\}^n$, we obtain
$$
0=\sum_{B\in\ka B} \sum_B f= \sum_{B\in [A]} \sum_B f= |[A]_e| \sum_Af-|[A]_o|\sum_Af,
$$ 
which completes the proof.\hfill $\square$

\medskip
Let us emphasize that we shall exploit only the second part of our lemma. 

For every $S\subseteq\{1,\ldots, n\}$, one defines the character $\chi_S\colon \{0,1\}^n\to \{-1,1\}$ by
$$
\chi_S(x)=(-1)^{\sum_{i\in S} x_i}. 
$$
Let us remark that the function $f$ defined in the course of the proof is simply equal to $\chi_{\prop A}$. (The reader is referred to \cite{O'Donnell} for further information on characters.) 

Every box $B\subseteq \{0,1\}^n$ can be encoded as a word $w=w_1w_2\cdots w_n$ over the alphabet $\{0, 1, *\}$ and conversely. The encoding is defined by the correspondence: $\{0\}\leftrightarrow 0$, $\{1\}\leftrightarrow 1$, $\{0,1\}\leftrightarrow *$. From now on, we shall use the terminology of boxes and words interchangeable. All notions considered so far, as for example function $B\mapsto \prop B$, translate to words in an obvious manner.

\begin{tw}
\label{glowne}
Let $3\le k<n$ be two integers. Let $V$ be a set of words of length $n$ over the alphabet $\{0,1,*\}$. Suppose the following conditions are satisfied:
\begin{itemize}
\item[$(\alpha_1)$]
$|v|=2^{n-k}$ for every $v\in V$ (equivalently, $|\prop v|=k$);
\item[$(\alpha_2)$]
if $v,u\in V$ are distinct, then there is exactly one $i$ such that $\{v_i,u_i\}=\{0,1\}$;
\item[$(\alpha_3)$]
if $v,u\in V$ are distinct, then $\prop u\neq \prop v$.
\end{itemize}
Then $|V| \le 2^k-2$.
\end{tw}

Two cases $k=1, 2$ are excluded from our theorem. The first of them is obvious: If $k=1$, then $|V|\le 1$. The following example shows that if $k=2$, then the upper bound for $|V|$ has to be at least $3$:
$$
\begin{array}{ccc}
\label{K2}
0&0&*\\
*&1&0\\
1&*&1
\end{array}\,.
$$ 
We shall show that it is 3. 

Let us start with elementary operations over words. We consider two types of such operations: those induced by permutations, and those induced by complementations: 
\begin{itemize}
\item[$(\alpha)$] If $\sigma$ is a permutation of the set $\{1,\ldots ,n\}$, then the operation over words of length $n$ induced by $\sigma$ is defined by $v\mapsto v\sigma =v_{\sigma(1)}\cdots v_{\sigma(n)}$. 

\item[$(\beta)$] Let $c\colon \{0,1,*\}\to \{ 0,1,*\}$ be given by $c(0)=1, c(1)=0, c(*)=*$. Every sequence $\gamma_1,\ldots,\gamma_n$, where each $\gamma_i$ is equal to $c$ or is the identity mapping on $\{0,1,*\}$ induces the mapping $v\mapsto\gamma(v)=\gamma_1(v_1)\cdots\gamma_n(v_n)$ defined on words of length $n$ over the alphabet $\{0,1,*\}$. 
\end{itemize}

It is clear that if $V$ is a set of words which fulfiles conditions $(\alpha_1$--$\alpha_3)$ of our theorem, then any set $V'$ which results from $V$ by consecutive applications of elementary operations also fulfiles $(\alpha_1$--$\alpha_3)$. The cardinality of $V'$ is equal to that of $V$. Therefore, we can always consider $V'$ instead of $V$ when we are looking for an estimate of $|V|$. We shall use such a replacement without further comments. 

Let us go back to the case $k=2$. We may assume without loss of generality that $u=00*\cdots*$ belongs to $V$ . By our assumptions, if $v\in V$ and $v\neq u$, then $*\in \{v_1,v_2\}$. We may assume that $v_1=*$. Then $v_2$ has to be $1$, as $v$ and $u$ has to fulfil $(\alpha_2)$. Moreover, we deduce from $(\alpha_1)$ that there is exactly one $i>2$ for which $v_i\in \{0,1\}$. Therefore, we may assume that $v=*10*\cdots*$. Now, it is easily seen that if $w\in V$ is distinct from $u$ and $v$, then our assumptions enforce $w$ to be equal to $1*1*\cdots*$.

\medskip
\noindent
\textit{Proof of the theorem.}

Let $\varepsilon \in \{0,1,*\}$. Let $V^{i\varepsilon}=\{ v\in V\colon v_i=\varepsilon\}$. 


\medskip
\noindent
\textbf{Claim 1.} \textit{If there is $i$ such that $|V^{i0}| \neq |V^{i1}|$, then $|V|\le 2^k-2$.}

We may assume that $i=1$ and $|V^{i0}|<|V^{i1}|$. Let us consider two words of length $n$: $\varepsilon*\cdots*$, $\varepsilon=0,1$. Let 
$$
W^{\varepsilon}=\{\varepsilon*\cdots* \cap v\colon \text{$v\in V$ and $\varepsilon*\cdots* \cap v\neq \emptyset$}\}.
$$
It is easily seen that if $x\in W^\varepsilon$ is minimal (with respect to $\preceq$), then , by our assumptions, $ W^\varepsilon$ does not contain any other element equivalent to $x$. Thus, by Lemma \ref{minbox}, boxes belonging to $W^\varepsilon$ cannot cover $\varepsilon*\cdots*$. Since the minimal cardinality of arbitrary box belonging to $W^\varepsilon$ is at least $2^{n-k-1}$, it follows that the uncovered part of $\varepsilon*\cdots*$ is a multiple of that number. The inequality $|V^{i0}|<|V^{i1}|$ implies that the uncovered part of $0*\cdots*$ is of greater cardinality than that of $1*\cdots*$. Thus, the uncovered by $V$ part of the box $*\cdots*$ is greater than $2^{n-k}$ and is a multiple of $2^{n-k}$. Consequently, it is at least $2^{n-k+1}$, which readily completes the proof of our claim. 

Therefore, we can further assume that $|V^{i0}|=|V^{i1}|$ for every $i$. Suppose now that for some $i$, one has $V^{i*}=\emptyset$. Since $V^{i\varepsilon}$ cannot cover $\varepsilon *\cdots*$, for $\varepsilon=0,1$, it appears that $|V^{i\varepsilon}|\le 2^{k-1}-1$. Then 
$$
|V|=|V^{i0}|+|V^{i1}|\le 2^k-2.
$$
Summing up, we may assume that 
$$
\label{addas}
\leqno{(A_1)}\qquad |V^{i0}|=|V^{i1}| \neq 0 \,\,\, \text{and}\,\,\,  V^{i*}\neq\emptyset, \,\,\, \text{for every $i\in \{1,\ldots,n\}$.} 
$$

We may also assume that $u=0\cdots 0*\cdots*\in V$. Clearly, $\prop u=\{1,\ldots, k\}$. Let $\delta$ be an arbitrary word of length $n-k$ over the alphabet $\{0,1\}$. Then $u^\delta=0\cdots 0\delta$ is a sub-box of $u$. (In fact, it is a singleton of an element of $u$). Consider a new word $*\cdots*\delta$ of length $n$. Let 
$$
A^\delta=\{v\in V\colon  v\cap*\cdots*\delta \neq \emptyset \},
$$ 
$$
B^\delta=\{v\cap*\cdots*\delta\colon v\in A^\delta \}.
$$
Since $u^\delta$ is an element of $B^\delta$, both sets $A^\delta$, $B^\delta$ are nonempty. Moreover, $u^\delta$ is a minimal (with respect to $\preceq$) element of $B^\delta$, and there is no other members of the latter set equivalent to $u^\delta$. By Lemma \ref{minbox}, there is an element $w^\delta\subseteq *\cdots*\delta$ which is disjoint with all members of $B^\delta$ so does with $V$, has an odd number $p_\delta$ of occurrences of `1' in first $k$ positions and is equivalent to $u^\delta$. Let 
$$
U^\delta=*\cdots*\delta\setminus \bigcup B^\delta.
$$
The set $U^\delta$ is the uncovered part of the $*\cdots*\delta$. Therefore $U=\bigcup_\delta U^\delta$, where the union extends over all words $\delta$ of length $n-k$ over the alphabet $\{0,1\}$, is the uncovered by $V$ part of the $n$-box $\{0,1\}^n=*\cdots*$. Since $w^\delta\subseteq U^\delta$, we have $|U^\delta|\ge 1$. Clearly, the sets $U^\delta$ are pairwise disjoint. Therefore, $|U|\ge 2^{n-k}$. We have to show that $|U|>2^{n-k}$ in order to complete the proof of our theorem. (In other words, we have to find a word $\tau$ of length $n-k$ over the alphabet $\{0,1\}$ so that $|U^\tau|>1$). Conversely, suppose that $U^\delta$ is a singleton for every $\delta$. 
Since $p_\delta$ is odd for every $\delta$, we can split our reasoning into two cases: (1) there is $\delta$ for which $p_\delta\ge3$; (2) $p_\delta=1$ for every $\delta$. 

\medskip\noindent
\textbf{Case 1.} We may assume that the first three symbols of word $w^\delta$ are `1's. Let us define three words of length $n$: $u^1=0*\cdots*\delta$, $u^2=*0*\cdots*\delta$ and $u^3=**0*\cdots*\delta$. Since $w^\delta$ is the only member of $U^\delta$, it follows that each of these words is disjoint with $U^\delta$. As $u^\delta=u\cap u^i$ for every $i$, by Lemma \ref{minbox}, there are elements $v^i\in V$ so that $v^i\cap u^i$ are different from but equivalent to $u^\delta$. Consequently, every $v^i$ has at most one star in the first $k$ positions, and if it occurs, then at position $i$. On the other hand, every $v^i$ must have this star, as otherwise they would be equivalent to$u$ which is forbidden by $(\alpha_3)$. It is also true that every of $v^i$ has exactly one symbol `1' in the first $k$ positions (It stems from $(\alpha_2)$ and the fact that $u$ has `0' in the first $k$ positions and stars in the remaining places). Moreover, if for a pair $v^i, v^j$ there is $s$ such that $\{v_s^i,v_s^j\}=\{0,1\}$, then necessarily $s\le k$, as the subwords $v^i_{k+1}\cdots v^i_{n}$, $v^j_{k+1}\cdots v^j_{n}$ contain $\delta$. Observe now that for every pair of different words $v^i$, $v^j$ they cannot have `1' at the same position $i \le k$, as if this is the case, then one of them should have two occurrences of `1' in the first $k$ positions, which is forbidden. It is also easily proved that `1' can occur only in the first three positions for every $v^i$ . Conversely, suppose $v^1$ has `1' at position 4, just to fix our attention. As $v^2$ and $v^3$ cannot have `1' at the same position, at least one of them, say $v^2$, has to have `1' at position $s$ different from the first and the fourth as well. Then $v^2$ and $v^1$ would violate $(\alpha_2)$, as $\{v^1_t,v^2_t\}=\{0,1\}$ for $t=s, 4$. Therefore, `1's can be distributed in one of the following two ways:
$$
\begin{array}{cc}
\begin{array}{ccc}
*&1&0\\
0&*&1\\
1&0&*
\end{array} &\qquad
\begin{array}{ccc}
*&0&1\\
1&*&0\\
0&1&*
\end{array}
\end{array}.
$$
As the reasoning is the same in both cases to be discussed, we shall consider only the first of them. Let $x$ be a word that belongs to $V^\star=V\setminus \{u,v^1,v^2, v^3\}$. We already know that $x$, as any other element of $V$ different from $u$, has to have a unique`1' in one of the first $k$ positions. Let us denote this position by $s$. If $s>3$, then $x$ begins with three stars in order to avoid violation of $(\alpha_2)$. If $s \le 3$, say for example $s=1$, then $x_2=*$, as otherwise $x_2=0$ and the pair $x$, $v^3$ would violate $(\alpha_2)$. Consider now pair $x$, $v^2$. There has to exist $s$ such that $\{x_s,v^2_s\}=\{0,1\}$. Clearly $s>k$. Now, from $(\alpha_1)$ and $(\alpha_3)$ we deduce that $x$ has to have an additional star in one of the first $k$ positions. Summing up, if $x\in V^\star$, then it has at least two stars in the first $k$ positions.

Let $p,q,r>k$ be these positions for which $v^1_p, v^2_q, v^3_r$ are different from `$*$'. Let us pick a word $\tau$ of length $n-k$ over the alphabet $\{0,1\}$ so that 
$$
\leqno{(A2)} \qquad \tau_{p-k}\neq v^1_p,\quad \tau_{k}\neq v^2_q,\quad  \tau_{r-k}\neq v^3_r.
$$
Let us consider the intersections of words belonging to $V$ with $*\cdots *\tau$, that is, the set $B^\tau$. Clearly, $u^\tau=u\cap *\cdots *\tau$ is a singleton. By (A2), 
none of the $v^i$ belongs to $A^\tau$. If $x\in V^\star$, then, by the fact that such an $x$ has at least two stars in the first $k$ positions, it follows that the cardinality of $x\cap*\cdots *\tau$ is a multiple of 4. Therefore, there is a unique element of $B^\tau$ which is of cardinality 1, while the others have their cardinalities divisible by 4. Since $*\cdots *\tau$ is a multiple of $8$, we conclude that $|U^\tau|$ is at least 3, which validates our theorem if the first case takes place. 

\medskip\noindent
\textbf{Case 2.}
Recall that every $w^\delta$ is a singleton of an element of $\{0,1\}^n$. Slightly abusing the terminology, we identify each $w^\delta$ with its only element. Then 
$$
U=\{w^\delta\colon \delta=\delta_1\cdots\delta_{n-k},\, \delta_1\in \{0,1\}, \ldots, \delta_{n-k}\in \{0,1\}\}.
$$ 
For $\varepsilon\in \{0,1\}$ and $i\in \{1,\ldots,k\}$, let $U^{i\varepsilon}=\{ w\in U\colon w_i=\varepsilon\}$. As is easily seen by $(A1)$, $|U^{i0}|=|U^{i1}|$ for every $i\in \{1,\ldots, k\}$. Since, by our assumption on $p_\delta$, symbol `1' appears only once in $w_1\cdots w_k$ for every $w \in U$, we conclude that $|U|= \sum_{i=1}^k |U^{i1}|$. On the other hand, `0' appears $k-1$ times in $w_1\cdots w_k$ for every $w \in U$, which shows that $(k-1)|U|= \sum_{i=1}^k |U^{i0}|$. Consequently, $|U|=(k-1)|U|$, which is impossible.
\hfill $\square$ 
\section{Conjecture}

Let $W$ be a set of words of a fixed length over the alphabet $\{0,1,*\}$. Suppose that $W$ satisfies $(\alpha_1$--$\alpha_3)$ with $k=m$. We already know that the maximum cardinality of $W$ if $k=2$ is 3. Let us define a new set of words
$$
W'=\{w*^n0\colon w \in W\}\cup\{*^nw1\colon w\in W\},
$$
where $*^n$ is the word consisting of $n$ stars `$*$'.
Clearly, $W'$ satisfies $(\alpha_1$--$\alpha_3)$ with $k=m+1$. Moreover, $|W'|=2|W|$. Therefore, we deduce by induction that for every $k\ge 2$ there is $W$ whose cardinality is $\frac{3}{4}2^k$. We conjecture that it is the maximum cardinality; that is, Theorem \ref{glowne} can be strengthen by replacing $2^k-2$ with $\frac{3}{4}2^k$. Let us remark that these two numbers coincide for $k=3$. Observe also that this conjecture implies $c_d\le \frac{3}{2}2^d$.

\skok
\noindent{\small
A.K.:  Wydzia{\l} Matematyki, Informatyki i Ekonometrii, Uniwersytet Zielonog\'orski, ul. Podg\'orna 50,\\ 65-246 Zielona G\'ora, Poland\\
{\tt A.Kisielewicz@wmie.uz.zgora.pl}
}
\skok
\noindent{\small
K.P.:  Wydzia{\l} Matematyki, Informatyki i Ekonometrii, Uniwersytet Zielonog\'orski, ul. Podg\'orna 50, \\
65-246 Zielona G\'ora, Poland\\
{\tt K.Przeslawski@wmie.uz.zgora.pl}}

\end{document}